\documentclass[11pt]{article}
\usepackage{latexsym,amssymb,amsmath,enumerate,float,geometry,cite,algorithm2e}
\newtheorem{theorem}{Theorem}
\newtheorem{corollary}[theorem]{Corollary}

\newtheorem{lemma}[theorem]{Lemma}
\newtheorem{conjecture}[theorem]{Conjecture}

\begin{document}

\title{\Large Bounds on the Burning Number}

\author{
St\'{e}phane Bessy$^1$\\ 
Anthony Bonato$^2$\\
Jeannette Janssen$^3$\\
Dieter Rautenbach$^4$\\
Elham Roshanbin$^3$
}
\date{}
\maketitle
\vspace{-10mm}
\begin{center}
{\small
$^1$ Laboratoire d'Informatique, de Robotique et de Micro\'{e}lectronique de Montpellier (LIRMM),\\ 
Montpellier, France, \texttt{stephane.bessy@lirmm.fr}\\[3mm]
$^2$ Department of Mathematics, Ryerson University, Toronto, ON,\\ 
Canada, M5B 2K3, \texttt{abonato@ryerson.ca}\\[3mm]
$^3$ Department of Mathematics and Statistics, Dalhousie University,\\
Halifax, NS, Canada, B3H 3J5, \texttt{jeannette.janssen,e.roshanbin@dal.ca}\\[3mm]
$^4$ Institute of Optimization and Operations Research, Ulm University,\\
Ulm, Germany, \texttt{dieter.rautenbach@uni-ulm.de}
}
\end{center}

\begin{abstract}
Motivated by a graph theoretic process intended to measure the speed of the spread of contagion in a graph, 
Bonato, Janssen, and Roshanbin 
[Burning a Graph as a Model of Social Contagion, Lecture Notes in Computer Science 8882 (2014) 13-22]
define the burning number $b(G)$ of a graph $G$
as the smallest integer $k$ for which there are vertices $x_1,\ldots,x_k$ such that 
for every vertex $u$ of $G$, there is some $i\in \{ 1,\ldots,k\}$ with ${\rm dist}_G(u,x_i)\leq k-i$,
and ${\rm dist}_G(x_i,x_j)\geq j-i$ for every $i,j\in \{ 1,\ldots,k\}$.

For a connected graph $G$ of order $n$, they prove that $b(G)\leq 2\left\lceil\sqrt{n}\right\rceil-1$, 
and conjecture $b(G)\leq \left\lceil\sqrt{n}\right\rceil$.
We show that
$b(G)\leq \sqrt{\frac{32}{19}\cdot \frac{n}{1-\epsilon}}+\sqrt{\frac{27}{19\epsilon}}$
and 
$b(G)\leq \sqrt{\frac{12n}{7}}+3\approx 1.309 \sqrt{n}+3$ 
for every connected graph $G$ of order $n$ and every $0<\epsilon<1$.
For a tree $T$ of order $n$ with $n_2$ vertices of degree $2$, and $n_{\geq 3}$ vertices of degree at least $3$, 
we show
$b(T)\leq \left\lceil\sqrt{(n+n_2)+\frac{1}{4}}+\frac{1}{2}\right\rceil$
and
$b(T)\leq \left\lceil\sqrt{n}\right\rceil+n_{\geq 3}$.
Furthermore, we characterize the binary trees of depth $r$ that have burning number $r+1$.
\end{abstract}

{\small
\begin{tabular}{lp{12.5cm}}
\textbf{Keywords:} & burning; distance domination\\
\textbf{MSC2010:} & 05C57; 05C69
\end{tabular}
}

\pagebreak

\section{Introduction}

Motivated by a graph theoretic process intended to measure the speed of the spread of contagion in a graph, 
Bonato, Janssen, and Roshanbin \cite{bjr1,bjr2} define a {\it burning sequence} of a graph $G$ 
as a sequence $(x_1,\ldots,x_k)$ of vertices of $G$ such that 
\begin{eqnarray}
\forall u\in V(G): \exists i\in [k]: & {\rm dist}_G(u,x_i)\leq & k-i\mbox{ and}\label{e1}\\
\forall i,j\in [k]: & {\rm dist}_G(x_i,x_j)\geq & j-i,\label{e2}
\end{eqnarray} 
where $[k]$ denotes the set of the positive integers at most $k$.
Furthermore, they define the {\it burning number} $b(G)$ of $G$ as the
length of a shortest burning sequence of $G$. 

A burning sequence is supposed to model the expansion of a fire within a graph: At each discrete time step, 
first a new fire starts at a vertex
that is not already burning, 
and then the fire spreads from burning vertices 
to all their neighbors that are not already burning. 
Condition (\ref{e1}) ensures
that putting fire to the vertices of a burning sequence 
$(x_1,\ldots, x_k)$ in the order $x_1,\ldots,x_k$, 
all vertices of $G$ are burning after $k$ steps.
Condition (\ref{e2}) ensures that one never puts fire 
to a vertex that is already burning.

We consider only finite, simple, and undirected graphs, and use
standard terminology and notation \cite{d}.  For a graph $G$, a vertex $u$ of
$G$, and an integer $k$, let $N_G^k[u]=\{ v\in V(G):{\rm
  dist}_G(u,v)\leq k\}$.  Note that $N_G^0[u]=\{ u\}$ and
$N_G^1[u]=N_G[u]=\{ u\}\cup N_G(u)$.

With this notation (\ref{e1}) is equivalent to 
\begin{eqnarray}
V(G) &=& N_G^{k-1}[x_1]\cup N_G^{k-2}[x_2]\cup \cdots \cup N_G^0[x_k].\label{e1b}
\end{eqnarray} 
As previously said, condition (\ref{e2}) is motivated by the
considered graph process, which in each step puts fire to a vertex 
that is not already burning. Our first result is that condition
(\ref{e2}) is redundant.

\begin{lemma}\label{lemma1}
The burning number of a graph $G$ is the minimum length of a sequence
$(x_1,\ldots,x_k)$ of vertices of $G$ satisfying (\ref{e1b}).
\end{lemma}
{\it Proof:} Let $k$ be the minimum length of a sequence satisfying (\ref{e1b}).
By definition, $b(G)\geq k$. It remains to show equality.
For a contradiction, suppose $b(G)>k$.
Let the sequence $s=(x_1,\ldots,x_k)$ be chosen such that (\ref{e1b}) holds,
and $j(s)=\min\{ j\in [k]:{\rm dist}_G(x_i,x_j)<j-i\mbox{ for some }i\in [j-1]\}$ is as large as possible.
Since $b(G)>k$, the index $j(s)$ is well defined.
Let $i(s)\in [j(s)-1]$ be such that ${\rm dist}_G(x_{i(s)},x_{j(s)})<j(s)-i(s)$.
Since $k>j(s)-1$, 
there is a vertex $y$ in $$V(G)\setminus \left(N_G^{(j(s)-1)-1}[x_1]\cup N_G^{(j(s)-1)-2}[x_2]\cup \cdots \cup N_G^{0}[x_{j(s)-1}]\right).$$
Since $N_G^{k-j(s)}[x_{j(s)}]\subseteq N_G^{k-i(s)}[x_{i(s)}]$,
the sequence $s'=(x_1,\ldots,x_{j(s)-1},y,x_{j(s)+1},\ldots,x_k)$ satisfies (\ref{e1b}) and $j(s')>j(s)$,
which is a contradiction.
$\Box$

\medskip

\noindent In view of Lemma \ref{lemma1}, 
the burning number can be considered a variation (but distinct from) of well known distance domination parameters \cite{h}.
For a graph $G$ and an integer $k$, a set $D$ of vertices of $G$ is a {\it distance-$k$-dominating set} of $G$
if $\bigcup_{x\in D}N_G^k[x]=V(G)$.
The {\it distance-$k$-domination number $\gamma_k(G)$} of $G$ is the minimum cardinality of a distance-$k$-dominating set of $G$.

The following bound on the distance-$k$-domination number will be of interest.
\begin{theorem}[Meir and Moon \cite{mm}]\label{theoremmm}
If $G$ is a connected graph of order $n$ at least $k+1$, 
then $\gamma_k(G)\leq\frac{n}{k+1}$.
\end{theorem}
As observed in \cite{bjr1,bjr2} the burning number can be bounded above in terms of the distance-$k$-domination number.
In fact, if $\{ x_1,\ldots,x_{\gamma}\}$ is a distance-$k$-dominating set of $G$, then 
\begin{eqnarray*}
V(G) &=& N_G^{k}[x_1]\cup N_G^{k}[x_2]\cup \cdots \cup N_G^{k}[x_{\gamma}]\\
& = & N_G^{k+\gamma-1}[x_1]\cup N_G^{k+\gamma-2}[x_2]\cup \cdots \cup N_G^{k}[x_{\gamma}].
\end{eqnarray*}
Appending any $k$ vertices to the sequence $(x_1,\ldots,x_{\gamma})$
yields a sequence of length $k+\gamma$ satisfying (\ref{e1b}),
which, by Lemma \ref{lemma1}, implies $b(G)\leq \gamma_k(G)+k$.
Using Theorem \ref{theoremmm} and choosing $k=\left\lceil\sqrt{n}\right\rceil-1$, 
this implies the following.

\begin{theorem}[Bonato, Janssen, and Roshanbin \cite{bjr1,bjr2}]\label{theoremjbr}
If $G$ is a connected graph of order $n$, 
then $b(G)\leq 2\left\lceil\sqrt{n}\right\rceil-1$.
\end{theorem}
One of the most interesting open problems concerning the burning number is the following.

\begin{conjecture}[Bonato, Janssen, and Roshanbin \cite{bjr1,bjr2}]\label{conjecturebjr}
If $G$ is a connected graph of order $n$, then $b(G)\leq \left\lceil\sqrt{n}\right\rceil$.
\end{conjecture}
Since the path $P_n$ of order $n$ has burning number $\left\lceil\sqrt{n}\right\rceil$ \cite{bjr1,bjr2},
the bound in Conjecture \ref{conjecturebjr} would be tight.

Let ${\rm rad}(G)$ denote the radius of a graph $G$. 
Since $V(G)=N_G^{{\rm rad}(G)}[x]$ for every connected graph $G$
and every vertex $x$ of $G$ of minimum eccentricity, 
Lemma \ref{lemma1} implies the following.

\begin{theorem}[Bonato, Janssen, and Roshanbin \cite{bjr2}]\label{theoremradius}
If $G$ is a connected graph, then $b(G)\leq {\rm rad}(G)+1$.
\end{theorem}
In the present note, we improve the bound of Theorem \ref{theoremjbr} by showing several upper bounds on the burning number,
thereby contributing to Conjecture \ref{conjecturebjr}.
Furthermore, we characterize the extremal binary trees for Theorem \ref{theoremradius}.

\section{Results}

We begin with two straightforward results that lead to a first improvement of Theorem \ref{theoremjbr},
and rely on arguments that are typically used to prove Theorem \ref{theoremmm}.
For a vertex $u$ of a rooted tree $T$, let $T_u$ denote the subtree of $T$ rooted in $u$ that contains $u$ as well as all descendants of $u$.
Recall that the height of $T_u$ is the eccentricity of $u$ in $T_u$.

\begin{lemma}\label{lemma2}
Let $T$ be a tree.
If the non-negative integer $d$ is such that $N_T^d[u]\not= V(T)$ for every vertex $u$ of $T$, 
then there is a vertex $x$ of $T$ and a subtree $T'$ of $T$ with $n(T')\leq n(T)-(d+1)$ and $V(T)\setminus V(T')\subseteq N_T^d[x]$.
\end{lemma}
{\it Proof:} Root $T$ at a vertex $r$. Since $N_T^d[r]\not= V(T)$, the height of $T$ is at least $d+1$.
The desired properties follow for a vertex $x$ such that $T_x$ has height exactly $d$ and the tree $T'=T-V(T_x)$.
$\Box$

\begin{theorem}\label{theorem1}
Let $T$ be a tree.
If the non-negative integers $d_1,\ldots,d_k$ are such that $\sum\limits_{i=1}^k(d_i+1)\geq n(T)$,
then there are vertices $x_1,\ldots,x_k$ of $T$ such that $\bigcup\limits_{i=1}^kN_T^{d_i}[x_i]=V(T)$.
\end{theorem}
{\it Proof:} For a contradiction, suppose that such vertices do no exist.
Repeatedly applying Lemma \ref{lemma2}, 
yields a sequence $x_1,\ldots,x_k$ of vertices of $T$
as well as a sequence $T_1,\ldots,T_k$ of subtrees of $T$ such that
$n(T_i)\leq n(T_{i-1})-(d_i+1)$ and $V(T_{i-1})\setminus V(T_i)\subseteq N_{T_{i-1}}^{d_i}[x_i]\subseteq N_T^{d_i}[x_i]$
for every $i\in [k]$, where $T_0=T$.
Note that after $j-1<k$ applications of Lemma \ref{lemma2},
our assumption implies that $N_{T_{j-1}}^{d_j}[u]\not= V(T_{j-1})$ for every vertex $u$ of $T_{j-1}$,
because otherwise
\begin{eqnarray*}
V(T) & \subseteq &
(V(T_0)\setminus V(T_1))\cup(V(T_1)\setminus V(T_2))\cup \cdots \cup (V(T_{j-2})\setminus V(T_{j-1}))\cup V(T_{j-1})\\
& \subseteq &   
\bigcup\limits_{i=1}^{j-1}N_T^{d_i}[x_i]\cup N_{T_{j-1}}^{d_j}[u]\\
& \subseteq & 
\bigcup\limits_{i=1}^{j-1}N_T^{d_i}[x_i]\cup N_{T}^{d_j}[u]
\end{eqnarray*}
for some vertex $u$ of $T$,
contradicting our assumption.
Therefore, the hypothesis of Lemma \ref{lemma2} remains satisfied throughout its repeated applications.
Now, $V(T)\setminus V(T_k)\subseteq \bigcup_{i=1}^kN_T^{d_i}[x_i]$.
Since $n(T_k)\leq n(T)-\sum_{i=1}^k(d_i+1)\leq 0$, it follows that $V(T_k)$ is empty,
again contradicting our assumption.
$\Box$

\medskip

\noindent The previous result already allows to improve Theorem \ref{theoremjbr}.
\begin{corollary}\label{corollary1}
If $G$ is a connected graph of order $n$, 
then $b(G)\leq \left\lceil \sqrt{2n+\frac{1}{4}}-\frac{1}{2}\right\rceil$.
\end{corollary}
{\it Proof:} If $H$ is a spanning subgraph of $G$, then $b(G)\geq b(H)$. 
Hence, we may assume that $G$ is a tree.
If $k=\left\lceil \sqrt{2n+\frac{1}{4}}-\frac{1}{2}\right\rceil$, 
then $((k-1)+1)+((k-2)+1)+\cdots+(0+1)={k+1\choose 2}\geq n(G)$.
By Theorem \ref{theorem1}, 
there are vertices $x_1,\ldots,x_k$ in $G$ with $\bigcup\limits_{i=1}^kN_G^{k-i}[x_i]=V(G)$.
By Lemma \ref{lemma1}, $b(G)\leq k$. $\Box$

\medskip

\noindent Note that Theorem \ref{theoremmm}
is tight for any graph that arises by attaching a path of order $k$ to each vertex of a connected graph. 
In fact, also Theorem \ref{theorem1} is tight for the same kind of graph.
Therefore, in order to further improve Theorem \ref{theoremjbr}, 
one really has to leverage the full spectrum of different distances associated with the different vertices in a burning sequence.
The following lemma offers some way of doing this.

\begin{lemma}\label{lemma3}
Let $T$ be a tree.
If the positive integers $d_1$ and $d_2$ are such that 
$d_2\geq \left\lceil\frac{3d_1}{2}\right\rceil$ and 
$N_T^{d_1}[u]\cup N_T^{d_2}[v]\not= V(T)$ for every two vertices $u$ and $v$ of $T$, 
then there are two vertices $x$ and $z$ of $T$ and a subtree $T'$ of $T$ 
with $n(T')\leq n(T)-\left(\left\lceil\frac{3d_1}{2}\right\rceil+d_2+2\right)$ and 
$V(T)\setminus V(T')\subseteq N_T^{d_1}[x]\cup N_T^{d_2}[z]$.
\end{lemma}
{\it Proof:} Root $T$ at a vertex $r$. 
Since $N_G^{d_2}[r]\not= V(T)$, the height of $T$ is at least $d_2+1$.
Let the vertex $z$ be such that $T_z$ has height exactly $d_2$.
Note that $V(T_z)\subseteq N_T^{d_2}[z]$ and 
$|V(T_z)|\geq d_2+1$.
Let $x$ be a descendant of $z$ such that ${\rm dist}_T(x,z)=d_2-d_1$
and $T_x$ has height exactly $d_1$.
Note that $d_2-d_1\geq \left\lceil\frac{d_1}{2}\right\rceil$.
Let the vertex $y$ on the path in $T$ between $x$ and $z$ be such that ${\rm dist}_T(x,y)=\left\lceil\frac{d_1}{2}\right\rceil$.

If $V(T_y)\subseteq N_T^{d_1}[x]$, 
then Lemma \ref{lemma2} applied to the tree $\tilde{T}=T-V(T_y)$ and the value $d_2$
implies the existence of a vertex $z'$ and a subtree $T'$ of $\tilde{T}$ 
with $n(T')\leq n(\tilde{T})-(d_2+1)$ and $V(\tilde{T})\setminus V(T')\subseteq N_{\tilde{T}}^{d_2}[z']$.
Now, we have that
\begin{eqnarray*}
n(T') & \leq & n(\tilde{T})-(d_2+1)\\
& = & n(T)-|V(T_y)|-(d_2+1)\\
& \leq & n(T)-\left(\left\lceil\frac{3d_1}{2}\right\rceil+d_2+2\right)
\end{eqnarray*}
and 
\begin{eqnarray*}
V(T)\setminus V(T') & = & (V(T)\setminus V(\tilde{T}))\cup (V(\tilde{T})\setminus V(T'))\\
& \subseteq & V(T_y)\cup N_{\tilde{T}}^{d_2}[z']\\
& \subseteq & N_T^{d_1}[x]\cup N_T^{d_2}[z'].
\end{eqnarray*}
Hence, we may assume that $V(T_y)\not\subseteq N_T^{d_1}[x]$.
This implies the existence of a descendant $y'$ of $y$ that is not a descendant of $x$ and satisfies ${\rm dist}_T(x,y')>d_1$.
By the choice of $x$, $y$, and $z$, this implies $|V(T_z)|\geq d_2+1+\left\lceil\frac{d_1}{2}\right\rceil$.
Lemma \ref{lemma2} applied to the tree $\tilde{T}=T-V(T_z)$ and the value $d_1$
implies the existence of a vertex $x'$ and a subtree $T'$ of $\tilde{T}$ 
with $n(T')\leq n(\tilde{T})-(d_1+1)$ and $V(\tilde{T})\setminus V(T')\subseteq N_{\tilde{T}}^{d_1}[x']$.
Now, we have that
\begin{eqnarray*}
n(T') & \leq & n(\tilde{T})-(d_1+1)\\
& = & n(T)-|V(T_z)|-(d_1+1)\\
& \leq & n(T)-\left(\left\lceil\frac{3d_1}{2}\right\rceil+d_2+2\right)
\end{eqnarray*}
and 
\begin{eqnarray*}
V(T)\setminus V(T') & = & (V(T)\setminus V(\tilde{T}))\cup (V(\tilde{T})\setminus V(T'))\\
& \subseteq & V(T_z)\cup N_{\tilde{T}}^{d_1}[x']\\
& \subseteq & N_T^{d_1}[x'] \cup N_T^{d_2}[z],
\end{eqnarray*}
which completes the proof. $\Box$

\begin{theorem}\label{theorem2}
If $G$ is a connected graph and $0<\epsilon<1$, then $b(G)\leq \sqrt{\frac{32}{19}\cdot \frac{n(G)}{1-\epsilon}}+\sqrt{\frac{27}{19\epsilon}}$.
\end{theorem}
{\it Proof:} As in the proof of Corollary \ref{corollary1}, we may assume that $G$ is a tree $T$.

Let $\ell=\left\lceil\log_9\left(\frac{3}{19\epsilon}\right)\right\rceil$.
Note that 
$\left(1-\frac{3}{19}\cdot \left(\frac{1}{9}\right)^{\ell}\right)\geq 1-\epsilon$ and $3^{\ell}<\sqrt{\frac{27}{19\epsilon}}.$
Let $k$ be the smallest integer such that 
$(1-\epsilon)\cdot \frac{19k^2}{32}+(1-\epsilon)\cdot \frac{3k}{8}\geq n(T)$
and 
$k\equiv 0$ ( mod $3^{\ell}$).
Note that $$k\leq 
\left\lceil\sqrt{\frac{32}{19}\cdot \frac{n(T)}{1-\epsilon}+\left(\frac{6}{19}\right)^2}-\frac{6}{19}+3^{\ell}-1\right\rceil
\leq \sqrt{\frac{32}{19}\cdot \frac{n(T)}{1-\epsilon}}+\sqrt{\frac{27}{19\epsilon}}.$$
For a contradiction, suppose that $b(G)>k$.

For $j\in [\ell]$, let $I_j=\left[\frac{2k}{3^j}-1\right]\setminus \left[\frac{k}{3^j}-1\right]=\left\{\frac{k}{3^j},\frac{k}{3^j}+1,\ldots,\frac{2k}{3^j}-1\right\}$.
Since $\frac{k}{3^j}$ is an integer, it follows that $\left\lceil\frac{3d}{2}\right\rceil\leq \frac{k}{3^j}+d$ for every $d\in I_j$.
Repeatedly applying Lemma \ref{lemma3} to the 
$\left(1-\frac{1}{3^{\ell}}\right)k$ disjoint pairs $\left\{d,\frac{k}{3^j}+d\right\}$ for $j\in [\ell]$ and $d\in I_j$,  
yields pairs of vertices $\left\{ x_d,x_{\frac{k}{3^j}+d}\right\}$ as well as a subtree $T'$ of $T$ such that
\begin{eqnarray*}
n(T') 
& \leq & n(T)
-\sum_{j=1}^{\ell}\sum_{d=\frac{k}{3^j}}^{\frac{2k}{3^j}-1}\left(\left\lceil\frac{3d}{2}\right\rceil+\left(\frac{k}{3^j}+d\right)+2\right)\\
& \leq & n(T)
-\sum_{j=1}^{\ell}\sum_{d=\frac{k}{3^j}}^{\frac{2k}{3^j}-1}\left(\frac{5d}{2}+\frac{k}{3^j}+2\right)\\
& = & n(T)
-\sum_{j=1}^{\ell}\left(\frac{1}{9^{j-1}}\cdot\frac{19k^2}{36}+\frac{1}{3^{j-1}}\cdot \frac{k}{4}\right)\\
& = & n(T)
-\left(1-\left(\frac{1}{9}\right)^{\ell}\right)\cdot \frac{19k^2}{32}
-\left(1-\left(\frac{1}{3}\right)^{\ell}\right)\cdot \frac{3k}{8}
\end{eqnarray*}
and
\begin{eqnarray*}
V(T)\setminus V(T') 
& \subseteq & \bigcup_{j=1}^{\ell}\bigcup_{d=\frac{k}{3^j}}^{\frac{2k}{3^j}-1}
\left(N_T^d[x_d]\cup N_T^{\left(\frac{k}{3^j}+d\right)}\left[x_{\frac{k}{3^j}+d}\right]\right)\\
& = & \bigcup_{i=\frac{k}{3^{\ell}}}^{k-1}N_T^i[x_i].
\end{eqnarray*}
Note that, similarly as in the proof of Theorem \ref{theorem1},
the assumption $b(G)>k$ implies that the hypothesis of Lemma \ref{lemma3} 
remains satisfied throughout its repeated applications.

Now, repeatedly applying Lemma \ref{lemma2} for all $\frac{k}{3^{\ell}}$ values $d$ in $\{ 0\}\cup [\frac{k}{3^{\ell}}-1]$,
yields vertices $x_0,\ldots, x_{\frac{k}{3^{\ell}}-1}$ and a subtree $T''$ of $T'$ such that
\begin{eqnarray*}
n(T'') & \leq & n(T')-\sum_{d=0}^{\frac{k}{3^{\ell}}-1}(d+1)\\
& = & n(T')-\left(\frac{1}{9}\right)^{\ell}\cdot\frac{k^2}{2}-\left(\frac{1}{3}\right)^{\ell}\cdot\frac{k}{2}
\end{eqnarray*}
and
$$V(T')\setminus V(T'')\subseteq\bigcup_{i=0}^{\frac{k}{3^{\ell}}-1}N_T^i[x_i].$$
Altogether, the vertices $x_0,\ldots,x_{k-1}$ satisfy 
$$V(T)\setminus V(T'')\subseteq\bigcup_{i=0}^{k-1}N_T^i[x_i].$$
Since 
\begin{eqnarray*}
n(T'') & \leq & n(T)
-\left(1-\left(\frac{1}{9}\right)^{\ell}\right)\cdot \frac{19k^2}{32}
-\left(1-\left(\frac{1}{3}\right)^{\ell}\right)\cdot \frac{3k}{8}
-\left(\frac{1}{9}\right)^{\ell}\cdot\frac{k^2}{2}
-\left(\frac{1}{3}\right)^{\ell}\cdot\frac{k}{2}\\
& = & n(T)
-\left(1-\frac{3}{19}\cdot \left(\frac{1}{9}\right)^{\ell}\right)\cdot \frac{19k^2}{32}
-\left(1+\left(\frac{1}{3}\right)^{\ell+1}\right)\cdot \frac{3k}{8}\\
& \leq & n(T)
-(1-\epsilon)\cdot \frac{19k^2}{32}
-(1-\epsilon)\cdot \frac{3k}{8}\\
& \leq & 0,
\end{eqnarray*}
it follows that $V(T'')$ is empty, which implies the contradiction $b(T)\leq k$. $\Box$

\medskip

\noindent Choosing in the above proof $\ell=1$, and $k$ as the smallest multiple of $3$ that satisfies 
$\frac{7}{12}k^2+\frac{5}{12}k\geq n(T)$,
allows to deduce a similar contradiction, which implies
$b(G)\leq \sqrt{\frac{12n(G)}{7}}+3\approx 1.309 \sqrt{n(G)}+3$ for every connected graph $G$.

The following results generalize the equality $b(P_n)=\left\lceil\sqrt{n}\right\rceil$,
and establish approximate versions of Conjecture \ref{conjecturebjr} under additional restrictions.

\begin{lemma}\label{lemma4}
If $n_1,\ldots,n_p$ and $k$ are positive integers such that
$n_1+\cdots+n_p+k(p-1)\leq k^2$, then $b(P_{n_1}\cup\cdots\cup
P_{n_p})\leq k$.
\end{lemma}
{\it Proof:} The proof is by induction on $n=n_1+\cdots+n_p$.  Let
$G=P_{n_1}\cup\cdots\cup P_{n_p}$ and $n_1\leq \ldots \leq n_p$.  Note
that $p\leq k$.  If $n_p\leq k-p+1$, let the set
$\{ x_1,\ldots,x_p\}$ contain a vertex from each component of $G$.
We have $V(G)=N_G^{k-1}[x_1]\cup N_G^{k-2}[x_2]\cup \cdots \cup
N_G^{k-p}[x_p]$, and Lemma \ref{lemma1} implies $b(G)\leq k$.  Hence,
we may assume that $n_p\geq k-p+2$, which implies $n\geq
(p-1)+(k-p+2)=k+1$.

If $n_p\geq 2k$, let $x_1$ be a vertex at distance $k-1$ from an
endvertex of a component of $G$ of order $n_p$. 
The graph
$G'=G-N^{k-1}_G[x_1]$ has $p$ components and $|N^{k-1}_G[x_1]|=2k-1$.
Since 
\begin{eqnarray*}
n_1+\cdots+n_{p-1}+(n_p-(2k-1))+(k-1)(p-1)
& \leq & n_1+\cdots+n_{p-1}+(n_p-(2k-1))+k(p-1)\\
&\leq & k^2-(2k-1)\\
&=& (k-1)^2,
\end{eqnarray*}
there are, by induction, vertices $x_2,\ldots,x_k$ such that 
\begin{eqnarray}\label{e3}
V(G')=N_{G'}^{(k-1)-1}[x_2]\cup N_{G'}^{(k-1)-2}[x_3]\cup \cdots \cup N_{G'}^0[x_k].
\end{eqnarray}
This implies (\ref{e1b}), and Lemma \ref{lemma1} implies $b(G)\leq k$.

Hence, we may assume that $n_p\leq 2k-1$.  In this case we choose
as $x_1$ a vertex of minimal eccentricity in a component of $G$ of
order $n_p$.  This implies that $G'=G-N^{k-1}_G[x_1]$ has $p-1$
components.  Since
\begin{eqnarray*}
n_1+\cdots+n_{p-1}+(k-1)(p-2) &\leq& k^2-n_p-\left(k(p-1)-(k-1)(p-2)\right)\\
& \leq & k^2-(k-p+2)-(k+p-2)\\
&=& k^2-2k\\
&<& (k-1)^2,
\end{eqnarray*}
there are, by induction, vertices $x_2,\ldots,x_k$ that satisfy (\ref{e3}),
which again implies $b(G)\leq k$. $\Box$

\medskip

\noindent Since $n+\left(\left\lceil\sqrt{n}\right\rceil+(p-1)\right)(p-1)\leq \left(\left\lceil\sqrt{n}\right\rceil+(p-1)\right)^2$
for positive integers $n$ and $p$, Lemma \ref{lemma4} implies the following.

\begin{corollary}[Roshanbin \cite{r}]\label{corroshanbin}
If the forest $T$ of order $n$ is the union of $p$ paths, then $b(T)\leq \left\lceil\sqrt{n}\right\rceil+(p-1)$.
\end{corollary}
We derive further consequences of Lemma \ref{lemma4}.

\begin{theorem}\label{theorem3}
If $T$ is a tree of order $n$ that has $n_{\geq 3}$ vertices of degree at least $3$, then $b(T)\leq \left\lceil\sqrt{n}\right\rceil+n_{\geq 3}$.
\end{theorem}
{\it Proof:} Clearly, we may assume that $n_{\geq 3}\geq 1$.
Let $k=\left\lceil\sqrt{n}\right\rceil+n_{\geq 3}$.
Let $x_1,\ldots,x_{n_{\geq 3}}$ be the vertices of degree at least $3$.
Let $T'=T-\{ x_1,\ldots,x_{n_{\geq 3}}\}$, and let 
$T''=T-N_T^{k-1}[x_1]\cup \cdots \cup N_T^{k-n_{\geq 3}}[x_{n_{\geq 3}}]$.
Every component of $T'$ is a path $P$ such that at least one endvertex of $P$ has a neighbor in $\{ x_1,\ldots,x_{n_{\geq 3}}\}$.
Therefore, the distinct components of $T''$ arise by removing at least $k-n_{\geq 3}=\left\lceil\sqrt{n}\right\rceil$ vertices 
from distinct components of $T'$.
This implies that if $T''=P_{n_1}\cup\cdots\cup P_{n_p}$, then 
$$n_1+\cdots+n_p+\left\lceil\sqrt{n}\right\rceil (p-1)<
\left(n_1+\left\lceil\sqrt{n}\right\rceil\right)+\cdots+\left(n_p+\left\lceil\sqrt{n}\right\rceil\right)
\leq n-n_{\geq 3}<\left\lceil\sqrt{n}\right\rceil^2.$$
Now, Lemma \ref{lemma4} implies the existence of vertices $y_1,\ldots,y_{\left\lceil\sqrt{n}\right\rceil}$ such that 
$$V(T'')=N_{T''}^{\left\lceil\sqrt{n}\right\rceil-1}[y_1]\cup \cdots \cup N_{T''}^{0}[y_{\left\lceil\sqrt{n}\right\rceil}].$$
We obtain
$$V(T)=N_T^{k-1}[x_1]\cup \cdots \cup N_T^{\left\lceil\sqrt{n}\right\rceil}[x_{n_{\geq 3}}]
\cup N_T^{\left\lceil\sqrt{n}\right\rceil-1}[y_1]\cup \cdots \cup N_T^{0}[y_{\left\lceil\sqrt{n}\right\rceil}],$$
and Lemma \ref{lemma1} implies $b(T)\leq k$. $\Box$

\begin{theorem}\label{theorem4}
If $T$ is a tree of order $n$ that has $n_2$ vertices of degree $2$, then $$b(T)\leq \left\lceil\sqrt{(n+n_2)+\frac{1}{4}}+\frac{1}{2}\right\rceil.$$
\end{theorem}
{\it Proof:} Let $k=\left\lceil\sqrt{(n+n_2)+\frac{1}{4}}+\frac{1}{2}\right\rceil$.
Note that $k(k-1)\geq n+n_2$. 
For a contradiction, suppose that $b(T)>k$.
Root $T$ at a vertex $r$.
As before, we obtain that the height of $T$ is at least $k$.
Let $x_d$ be a vertex of $T$ such that $T_{x_d}$ has height exactly $d$ for some $d\in \{ 0\}\cup [k-1]$.
Let $V(T_{x_d})$ contain exactly $p_d$ vertices that have degree $2$ in $T$.
If $P$ is a path of length $d$ between $x_d$ and a leaf of $T_{x_d}$, 
then at least $d-p_d$ vertices of $P$ have a child that does not lie on $P$.
Therefore, $|V(T_{x_d})\setminus \{ x_d\}|\geq 2d-p_d$, and 
$T'=T-(V(T_{x_d})\setminus \{ x_d\})$ is a tree with $n_2-p_d$ vertices of degree $2$
such that $V(T)\setminus V(T')\subseteq N_T^d[x_d]$.
Note that $x_d$ has degree $1$ in $T'$.
Iteratively repeating this argument similarly as in the previous proofs, 
we obtain vertices $x_0,\ldots,x_{k-1}$ and integers $p_0,\ldots,p_{k-1}$ 
such that $p_0+\cdots+p_{k-1}\leq n_2$ and $\sum_{d=0}^{k-1}(2d-p_d)\leq n$.
Since $\sum_{d=0}^{k-1}(2d-p_d)\geq k(k-1)-n_2\geq n$,
we obtain $V(T)=N_T^0[x_0]\cup\cdots\cup N_T^{k-1}[x_{k-1}]$,
which implies the contradiction $b(G)\leq k$. $\Box$ 

\medskip

\noindent In view of the simple argument that shows Theorem \ref{theoremradius},
the extremal graphs for this bound might have a rather special structure.
Our final result supports this intuition for binary trees.

Recall that a rooted tree is {\it binary} if every vertex has at most two children,
and that a binary tree is {\it perfect} if every non-leaf vertex has exactly two children, and all leaves have the same depth,
that is, the same distance from the root.
Let $T_1$ be the rooted tree of order $2$, and,
for an integer $r$ at least $2$, 
let $T_r$ be the rooted tree that arises from the perfect binary tree of depth $r-1$ by subdividing all edges that are incident with a leaf.
Alternatively, $T_r$ arises by attaching a new leaf to each of the $2^{r-1}$ leaves of the perfect binary tree of depth $r-1$.

\begin{theorem}\label{theoremradiusextr}
If $r$ is a positive integer and $T$ is a binary tree of depth $r$, then $b(T)=r+1$ if and only if $T$ contains $T_r$ as a subtree.
\end{theorem}
{\it Proof:} Since the statement is trivial for $r=1$, we may assume that $r\geq 2$.

First, we show that $T=T_r$ has burning number $r+1$.
For a contradiction, suppose that $b(T)\leq r$.  
Let $u$ be the root of $T$, and let $v^1$ and $v^2$ be the two children of $u$.  
For $i\in [2]$, let $T^i$ be the subtree of $T$ rooted in $v^i$ that contains $v^i$ as well as all descendants of $v^i$ in $T$.
By Lemma \ref{lemma1}, there are vertices $x_1,x_2,\ldots,x_r$ with $V(T)=N_T^{r-1}[x_1]\cup N_T^{r-2}[x_2]\cup \cdots \cup N_T^0[x_r]$.
By symmetry, we may assume that $x_1\not\in V(T^1)$.  
Let $L$ be the set of leaves of $T$ that belong to $T^1$.  
Since $T^1$ is isomorphic to $T_{r-1}$, we have $|L|=2^{r-2}$.
Note that $N_G^{r-1}[x_1]$ does not intersect $L$.  
Furthermore, for every $i\in [r-1]\setminus \{ 1\}$, 
the set $N_T^{r-i}[x_i]$ contains at most $2^{r-i-1}$ vertices from $L$. 
In fact, the set $N_T^{r-i}[x_i]$ contains exactly $2^{r-i-1}$ vertices from $L$ if and only if $x_i\in V(T^1)$ and $x_i$ has depth $i$ in $T$. 
Since $N_T^0[x_r]=\{ x_r\}$, the set $N_T^0[x_r]$ contains at most one vertex from $L$.
Since $|L|=2^{r-2}=\sum_{i=2}^{r-1}2^{r-i-1}+1$, 
every vertex in $L$ belongs to exactly one of the sets $N_T^{r-i}[x_i]$ for $i=[r]\setminus \{ 1\}$.
This implies that $x_2,\ldots,x_r\in V(T^1)$, 
$x_i$ has depth $i$ in $T$ for $i\in [r-1]\setminus \{ 1\}$, and
$x_r$ is a leaf of $T$.
Let $u_0\ldots u_r$ be the path in $T$ from the root $u=u_0$ to the leaf $x_r=u_r$.
Note that $u_1=v^1$.
Since $x_2$ belongs to $T^1$, $x_2$ has depth $2$ in $T$, and $x_r\not\in N_T^{r-2}[x_2]$, 
the vertex $x_2$ is the child of $u_1$ distinct from $u_2$.
Moreover, as every vertex of $L$ belongs to exactly one of the sets $N_T^{r-i}[x_i]$ for $i\in [r]\setminus \{ 1\}$, 
no vertex $x_i$ with $i\in [r]\setminus \{ 1,2\}$ is a descendant of $x_2$. 
Iterating these arguments, it follows that, for every $i\in [r-1]\setminus \{ 1\}$, 
the vertex $x_i$ is the child of $u_{i-1}$ distinct from $u_i$. 
However, this implies the contradiction that $u_{r-1}\not\in N_T^{r-1}[x_1]\cup N_T^{r-2}[x_2]\cup \cdots \cup N_T^0[x_r]$.
Altogether, we obtain that $T_r$ has burning number $r+1$.
Together with Theorem \ref{theoremradius}, 
this implies that a binary tree $T$ of depth $r$ has burning number $r+1$ if $T$ contains $T_r$ as a subtree.

For the converse, we assume that $T$ is a binary tree of depth $r$ that does not contain $T_r$ as a subtree.
It follows 
that $T$ has a leaf of depth less than $r$ 
or 
that $T$ has a vertex of depth less than $r-1$ that has only one child.
In both cases we will show that $b(T)\leq r$.
First, we assume that $T$ has a leaf at depth less than $r$. 
Let $d$ be the minimum depth of a leaf of $T$. 
Let $u_0\ldots u_d$ be a path in $T$ between the root $u_0$ and a leaf $u_d$.  
By assumption, we have $d<r$.  
For $i\in [d]$, let $x_i$ be the child of $u_{i-1}$ that is distinct from $u_i$.  
Note that the subtree of $T$ rooted in $x_i$ that contains $x_i$ as well as all descendants of $x_i$ in $T$ has depth at most $r-i$.  
This implies that $V(T)=N_T^{r-1}[x_1]\cup N_T^{r-2}[x_2]\cup \cdots \cup N_T^{r-d}[x_d]\cup N_T^0[u_d]$, 
and, by Lemma \ref{lemma1}, we obtain $b(T)\leq r$.
Next, we assume that $T$ has a vertex $x$ of depth less than $r-1$ that has only one child.
Let $T'$ arise from $T$ by adding a new leaf $y$ as a child of $x$.
Clearly, $T'$ is a binary tree of depth $r$ that has a leaf of depth less than $r$, 
and, hence, $b(T)\leq b(T')\leq r$.
$\Box$

\medskip

\noindent {\bf Acknowledgment} This paper is part of a collaborative work that grew out of \cite{br} and \cite{bjr3}.

\end{document}